\documentclass[10]{amsart}

\usepackage{url}
\usepackage{latexsym}
\usepackage{lscape} 
\usepackage{pdflscape}
\usepackage{float}
\usepackage{mathtools}
\usepackage{amssymb}
\usepackage{xcolor}
\usepackage{diagrams}
\usepackage{tikz}
\newtheorem{thm}[equation]{Theorem}

\newtheorem{ex}[equation]{Example}
\newtheorem{pbm}[equation]{Problem}

\newtheorem{pro}[equation]{Proposition}

\theoremstyle{remark}
\newtheorem{rem}[equation]{Remark}

\theoremstyle{definition}
\newtheorem{Def}[equation]{Definition}

\numberwithin{equation}{section}

\newcommand{\Hom}{\text{Hom}}

\newcommand{\Ad}{\text{Ad}}

\newcommand{\fl}{{\mathfrak l}}

\newcommand{\fq}{{\mathfrak q}}

\newcommand{\fs}{{\mathfrak s}}

\newcommand{\be}{\begin{equation}}
\newcommand{\beu}{\begin{equation*}}

\newcommand{\spec}{{\mathrm{spec}}}

\begin{document}

\baselineskip=18pt
\sloppy
\title[Spectrum of locally symmetric spaces]
{Spectrum of semisimple locally symmetric spaces\\ and admissibility of spherical representations}

\author{S. Mehdi}
\author{M. Olbrich}
\address{Institut Elie Cartan de Lorraine, CNRS - UMR 7502, Universit\'e de Lorraine -  Campus de Metz, France}
\email{salah.mehdi@univ-lorraine.fr}
\address{Facult\'e des Sciences, de la Technologie et de la Communication , Universit\'e du Luxembourg, Luxembourg}
\email{martin.olbrich@uni.lu}
\maketitle

\begin{abstract}
{\it We consider compact locally symmetric spaces $\Gamma\backslash G/H$ where $G/H$ is a non-compact semisimple symmetric space and $\Gamma$ is a discrete subgroup of $G$. We discuss some features of the joint spectrum of the (commutative) algebra $D(G/H)$ of invariant differential operators acting, as unbounded operators, on the Hilbert space $L^2(\Gamma\backslash G/H)$ of square integrable complex functions on $\Gamma\backslash G/H$. In the case of the Lorentzian symmetric space $SO_0(2,2n)/SO_0(1,2n)$, the representation theoretic spectrum is described explicitly. The strategy is to consider connected reductive Lie groups $L$ acting transitively and co-compactly on $G/H$, a cocompact lattice $\Gamma\subset L$, and study the spectrum of the algebra $D(L/L\cap H)$ on $L^2(\Gamma\backslash L/L\cap H)$. Though the group $G$ does not act on $L^2(\Gamma\backslash G/H)$, we explain how (not necessarily unitary) $G$-representations enter into the spectral decomposition of $D(G/H)$ on $L^2(\Gamma\backslash G/H)$ and why one should expect a continuous contribution to the spectrum in some cases. As a byproduct, we obtain a result on the $L$-admissibility of $G$-representations. These notes contain the statements of the main results, the proofs and the details will appear elsewhere.}
\end{abstract}

\section{Introduction.}\label{Intro}

Let $G$ be a connected non-compact semisimple real Lie group and $H$ a connected closed subgroup of $G$ with complexified Lie algebras ${\mathfrak g}$ and ${\mathfrak h}$ respectively. Suppose that $G/H$ is a semisimple symmetric space of the non-compact type with respect to an involution $\sigma$ and consider the associated decomposition 
\begin{equation*}
{\mathfrak g}={\mathfrak h}\oplus{\mathfrak q}\;\text{ where }{\mathfrak q}:=\{X\in{\mathfrak g}\mid\sigma(X)=-X\}
\end{equation*}
Fix an invariant measure on $G/H$ and consider the Hilbert space $L^2(G/H)$ of square integrable complex functions on $G/H$. The action by left translations of $G$ on $L^2(G/H)$ defines a unitary representation whose decomposition is known as the Plancherel formula for symmetric spaces \cite{BS1}\cite{BS2}\cite{Del}\cite{MO}:
\begin{equation*}
L^2(G/H)\simeq \int^\oplus_{\widehat{G}_H}\big( W^H_{\widetilde{\rho},-\infty}\big)_0\otimes W_\rho d\mu(\rho)
\end{equation*}
\hspace*{8.3cm}$\circlearrowleft$\hspace*{1.4cm}$\circlearrowleft$
\vspace*{-0.1cm}\\
\hspace*{7.7cm}$D(G/H)$\hspace*{0.1cm}$\twoheadleftarrow$\hspace*{0.3cm}${\mathcal Z}({\mathfrak g})$
\vspace*{0.2cm}\\
where $W_\rho$ is a unitary irreducible representation of $G$, $W_{\widetilde{\rho}}$ the dual representation, $W_{\widetilde{\rho},-\infty}$ the space of distribution vectors of $W_{\widetilde{\rho}}$, $W^H_{\widetilde{\rho},-\infty}$ the space of $H$-invariants distribution vectors, $\big( W^H_{\widetilde{\rho},-\infty}\big)_0$ is a finite dimensional Hilbert space embedded in $W^H_{\widetilde{\rho},-\infty}$ whose dimension is the multiplicity of $W_\rho$ in $L^2(G/H)$, $\widehat{G}_H$ is the $H$-spherical dual $\{\rho\in\widehat{G}\mid W^H_{\widetilde{\rho},-\infty}\neq\{0\}\}$ and $d\mu$ the Plancherel measure. The (commutative) algebra $D(G/H)$ of $G$-invariant differential operators on $G/H$ acts on $\big(W^H_{\widetilde{\rho},-\infty}\big)_0$ and the center ${\mathcal Z}({\mathfrak g})$ of the enveloping algebra ${\mathcal U}({\mathfrak g})$ of ${\mathfrak g}$ acts via infinitesimal character on $W_\rho$. Except for exceptional cases where $G$ is $E_6$, $E_7$ and $E_8$ which we will not consider here, there is a surjective map ${\mathcal Z}({\mathfrak g})\twoheadrightarrow D(G/H)$ \cite{Hel}. In other words, harmonic analysis on $G/H$ is closely related to the spectrum of $D(G/H)$ on $L^2(G/H)$: 
\vspace*{0.1cm}
\begin{center}
harmonic analysis on $G/H$ $\longleftrightarrow$ $\text{Spec}_{L^2(G/H)}(D(G/H))$
\end{center}

A compact locally symmetric space is a smooth compact manifold of the form $\Gamma\backslash G/H$ where $\Gamma$ is a discrete subgroup of $G$, in particular the natural map $G\rightarrow \Gamma\backslash G/H$ is smooth for the quotient topology. As above one may ask about the analysis on 
$\Gamma\backslash G/H$. Let us consider the following group case example.
\begin{ex}
$G=G_1\times G_1$, $H=\Delta(G_1\times G_1)$ and $\Gamma=\Gamma_1\times\{e\}$, where $G_1$ is a connected non-compact semisimple real Lie group, $\sigma$ is the permutation $G_1\times G_1\rightarrow G_1\times G_1,\;(a,b)\mapsto (b,a)$ and $\Gamma_1$ a co-compact torsion free discrete subgroup of $G_1$. Then $L^2(G/H)\simeq L^2(G_1)$ is described by the Harish-Chandra Plancherel formula and 
\begin{equation*}
L^2(\Gamma\backslash G/H)\simeq L^2(\Gamma_1\backslash G_1)\simeq\bigoplus_{\pi\in\widehat{G}_1}V^{\Gamma_1}_{\widetilde{\pi},-\infty}\widehat{\otimes}V_\pi
\end{equation*}
where $\dim (V^{\Gamma_1}_{\widetilde{\pi},-\infty})<\infty$, $\dim (V_\pi)=\infty$ (except for the trivial representation $\pi$) and the set $\{\pi\in\widehat{G}_1\mid V^{\Gamma_1}_{\widetilde{\pi},-\infty}\neq\{0\}\}$ is discrete. Moreover $D(G/H)\simeq{\mathcal Z}({\mathfrak g})$ acts scalarly on $\pi$ by infinitesimal characters. In particular, the spectrum of $D(G/H)$ on $L^2(\Gamma\backslash G/H)$ is discrete, consists of eigencharacters only, and (almost) all eigenspaces are infinite dimensional.
\end{ex}
\noindent
Already this example shows that analysis on $\Gamma\backslash G/H$ is often related to automorphic forms:
\begin{center}
analysis on $\Gamma\backslash G/H$ $\longleftrightarrow$ automorphic forms
\end{center}
\begin{pbm}\label{thepbm}
Describe the spectrum $\text{Spec}_{L^2(\Gamma\backslash G/H)}(D(G/H))$ of $D(G/H)$ on $L^2(\Gamma\backslash G/H)$ in general.
\end{pbm}
{\bf Issues:}
\begin{itemize}\label{issues}
\item[(1)] When $H$ is not compact then a discrete subgroup $\Gamma\subset G$ such that $\Gamma\backslash G/H$ is a smooth compact manifold need not exist as it is illustrated by the Calabi-Markus phenomenon for $SO_0(1,n+1)/SO_0(1,n)$ \cite{CM}. Given such a $\Gamma$, the space $\Gamma\backslash G/H$ will be called a compact {\it Clifford-Klein form} of $G/H$.
One can ask for deformations of $\Gamma$ such that the new quotient is still a Clifford-Klein form (see \cite{Kos01}).
\item[(2)] The group $G$ does not act on $\Gamma\backslash G/H)$, and it is not clear how representations of $G$ are involved in $\text{Spec}_{L^2(\Gamma\backslash G/H)}(D(G/H))$.
\item[(3)] When $H$ is not compact, the coset $\Gamma\backslash G$ is not compact and it is not clear if $L^2(\Gamma\backslash G)$ and $L^2(\Gamma\backslash G/H)$ are related.
\item[(4)] $D(G/H)$ acts on $C^\infty(\Gamma\backslash G/H)$ and $C^\infty(\Gamma\backslash G/H)$ is dense in $L^2(\Gamma\backslash G/H)$, so $D(G/H)$ acts via unbounded operators on $L^2(\Gamma\backslash G/H)$. Some care is required with the domains of the operators. 
\end{itemize}

Problem \ref{thepbm} is related with a program developed by Kobayashi to study hidden symmetries (see \cite{Kos17}\cite{Kos09} and references therein). 
In the following we shall describe and state some of the main results we prove in \cite{MOl}.

\section{Existence of compact Clifford-Klein forms and spherical triples.}\label{spherical}

When $H$ is compact, the existence of a discrete subgroup $\Gamma$ of $G$ such that $\Gamma\backslash G/H$ is a smooth compact manifold is equivalent to the existence of torsion free uniform lattices in $G$. It turns out that torsion free uniform lattices always exist in $G$ by a result of Borel and Harish-Chandra \cite{BH}. 

For $H$ noncompact, this gives a sufficient condition for the existence of $\Gamma$, employed extensively by Kobayashi starting with \cite{Kos89}:
\begin{equation*}
\left.
    \begin{array}{lll}
      \text{existence of a closed connected
Lie group $L\subset G$ such that}\\
      (i) \text{ $L$ is reductively embedded in $G$}\\
      (ii) \text{ $L$ acts transitively on $G/H$}\\
      (iii) \text{ $L\cap H$ is compact}
    \end{array}
  \right\}  \Longrightarrow 
  \left.
    \begin{array}{l}
      \text{existence of compact Clifford-Klein}\\
        \text{forms for $G/H$}\\
    \end{array}
  \right.
\end{equation*}
Indeed, if $\Gamma$ is a co-compact lattice in $L$ then $\Gamma\backslash G/H$ is a compact locally symmetric space. Triples $(G,H,L)$ satisfying (i) and (ii) above were classified, when $G$ is simple, by Onishik in the late 1960's \cite{Oni69} (see also \cite{KY05}). We will refer to triples $(G,H,L)$ satisfying (i)(ii)(iii) 
as {\it transitive triples}. 

In fact, the above implication remains true, if we weaken {\em transitively} in (ii) to {\em cocompactly}. However, we are able to show  that this does not
change anything:  {\em cocompactly} implies  {\em transitively} \cite{MOl}.

Next fix $P=MAN$ the Langlands decomposition of a minimal parabolic subgroup in $G$ and consider a principal series representation $\text{Ind}_P^G \tau\otimes e^\nu\otimes 1$ of $G$, where $\tau\in\widehat{M}$ is a finite dimensional irreducible representation of $M$ and $\nu$ is a linear form on the complexified Lie algebra of $A$. Then van den Ban proved in \cite{Ban87} that the space of $H$-invariants distributions vectors in $\text{Ind}_P^G \tau\otimes e^\nu\otimes 1$ is finite dimensional:
\begin{equation*}
\dim\big(\text{Ind}_P^G \tau\otimes e^\nu\otimes 1\big)^H_{-\infty}<\infty
\end{equation*}
Now pick a transitive triple $(G,H,L)$ and let $P_L=M_LA_LN_L$ be the Langlands decomposition of a minimal parabolic in $L$. Let $\tau_L\in\widehat{M_L}$ be a finite dimensional irreducible representation of $M_L$ and $\nu_L$ a linear form on the complexified Lie algebra of $A_L$. Observe that $L$ is not 
$\sigma$-stable, $L/ L\cap H$ need not be a symmetric space and the space of $L\cap H$-invariants distributions $\big(\text{Ind}_{P_L}^L \tau_L\otimes e^{\nu_L}\otimes 1\big)^{L\cap H}_{-\infty}$ could be infinite dimensional as it is the case for the triple $(G_1\times G_1,\Delta(G_1\times G_1),G_1\times\{e\})$. In fact we obtain the following criterion for the space of $L\cap H$-invariants to be finite dimensional. 
\begin{pro}\label{inv} For any $\tau_L\in\widehat{M_L}$ a finite dimensional irreducible representation of $M_L$ and $\nu_L$ a linear form on the complexified Lie algebra of $A_L$, one has
$$\dim\big(\text{Ind}_{P_L}^L \tau_L\otimes e^{\nu_L}\otimes 1\big)^{L\cap H}_{-\infty}<\infty\Leftrightarrow P_L\text{ acts transitively on }L/L\cap H.$$
\end{pro}
\noindent
Transitive triples $(G,H,L)$ such that $P_L$ acts transitively on $L/L\cap H$ will be called {\it spherical triples}. The use of the qualifier `spherical' is motivated by the following facts:
\begin{eqnarray*}
L/L\cap H\text{ is a spherical homogeneous space}&\xLeftrightarrow[\hspace*{1.3cm}]{\;\;\;}&P_L\text{ acts with an open orbit in }L/L\cap H\\
&\xLeftrightarrow[\text{is compact}]{\text{if $L\cap H$}}&P_L\text{ acts transitively on }L/L\cap H
\end{eqnarray*}
%
\newpage
\begin{ex}\mbox{}\vspace*{-0.7cm}\label{obst}\\
\begin{itemize}
\item[(1)] $(G_1\times G_1,\Delta(G_1\times G_1), G_1\times\{e\})$ is not spherical.
\item[(2)] $(G_1\times G_1,\Delta(G_1\times G_1), G_1\times K_1)$ is spherical (where $K_1$ is maximal compact in $G_1$).
\item[(3)] $(SO_0(2,2n),SO_0(1,2n),U(1,n))$ is spherical ($n\geq 2$).
\item[(4)] $(SO_0(4,3),SO_0(4,1)\times SO(2),G_{2(2)})$ is not spherical.
\end{itemize}
\end{ex}
%

\section{$L$-admissibility.}\label{admissibility}
Suppose $G^\prime\subset G$ is a connected closed reductive subgroup and $\pi$ is an irreducible unitary representation of $G$. 
From a general result of Mautner and Teleman, the restriction of $\pi$ to $G^\prime$ can be decomposed as the direct integral sum of irreducible unitary representations: 
\begin{equation*}
\pi_{\mid_G{^\prime}}\simeq\int_{\widehat{G^\prime}}^\oplus M_\pi(\tau)\widehat{\otimes} V_\tau \;d\mu(\tau)
\end{equation*}
where $M_\pi(\tau)$ is the multiplicity (Hilbert) space.
\begin{Def}
$\pi$ is $G^\prime$-admissible if $\pi_{\mid_G{^\prime}}$ decomposes discretely and $\dim M_\pi(\tau)<\infty$
\end{Def}
\begin{ex}\mbox{}\vspace*{-0.7cm}\\
\begin{itemize}
\item[(1)] $\widehat{G}\ni\pi$ is $K$-admissible (Harish-Chandra admissibility theorem \cite{Har}).
\item[(2)] $G=SL(2,{\mathbb R})\times SL(2,{\mathbb R})$ and $H=\Delta(SL(2,{\mathbb R})\times SL(2,{\mathbb R}))$. 
Let $\widehat{G}_H\ni\pi=\pi_1\otimes\pi_1$ where $\pi_1$ is a holomorphic discrete series representation of $SL(2,{\mathbb R})$. From \cite{Rep}, 
$\pi_{\mid_{H}}$ has continuous spectrum (but multiplicities are finite). Therefore $\pi$ is not $H$-admissible.
\end{itemize}
\end{ex}
There is a criterion for $H$-admissibility in terms of associated varieties due to Kobayashi \cite{Kos98}. 
We prove the following admissibility result for spherical representations associated with triples $(G,H,L)$.
\begin{thm}\label{ladmiss}
Suppose $(G,H,L)$ is a spherical triple. If $\pi\in\widehat{G}_H$ then $\pi$ is $L$-admissible. Moreover, if $\pi$ is not trivial then there are infinitely many summands in $\pi_{\mid_L}$ all having $L\cap H$-invariants.
\end{thm}
%

\section{Embedding of Casimir operators.}\label{casimir}

Fix a transitive triple $(G,H,L)$. Write ${\mathcal U}({\mathfrak g})^H$ for the $H$-invariant elements in the enveloping algebra 
of ${\mathfrak g}$ and ${\mathcal U}({\mathfrak g}){\mathfrak h}$ the left ${\mathcal U}({\mathfrak g})$-ideal generated by ${\mathfrak h}$. The algebra $D(G/H)$ of $G$-invariant differential operators on $G/H$ is isomorphic to the following quotient algebra:
\begin{equation*}
D(G/H)\simeq {\mathcal U}({\mathfrak g})^H/{\mathcal U}({\mathfrak g}){\mathfrak h}\cap{\mathcal U}({\mathfrak g})^H
\end{equation*}
It is known that since $G/H$ is a symmetric space, $D(G/H)$ is commutative (see \cite{Hel00} and references therein). 

By definition, if $(G,H,L)$ is a transitive triple then one has a diffeomorphism of homogeneous spaces 
\begin{equation*}
G/H\simeq L/L\cap H
\end{equation*}
and an isomorphism of algebra 
\begin{equation*}
{\mathcal U}({\mathfrak g})\simeq{\mathcal U}({\mathfrak h})\otimes_{{\mathcal U}({\mathfrak l}\cap{\mathfrak h})}{\mathcal U}({\mathfrak l}).
\end{equation*}
In particular, we get an embedding of algebras
\begin{equation*}
\imath:D(G/H)\hookrightarrow D(L/L\cap H)\ .
\end{equation*}
The algebra $D(L/L\cap H)$ need not be commutative, though it is for many transitive triples (see \cite{KK19}). Of special interest in $D(G/H)$ (resp. $D(L/L\cap H)$, $D(L\cap K/L\cap H)$) is the Casimir element $\Omega_G$ (resp. $\Omega_L$, $\Omega_{L\cap K}$). The algebras $D(L/L\cap H)$ and $D(L\cap K/L\cap H)$ are defined in an analogous way to $D(G/H)$. Observe that when $G/H$ has rank one then $D(G/H)$ is generated by the Casimir $\Omega_G$. 
\begin{pro}
For each transitive triple $(G,H,L)$ with $G/H$ irreducible, the embedding $\imath(\Omega_G)$  is computed explicitly in terms of the decomposition of ${\mathfrak l}/{\mathfrak l}\cap{\mathfrak h}$ into $\Ad(L\cap H)$-irreducibles.
\end{pro}
\noindent
In particular, based on the restriction to ${\mathfrak l}$ and ${\mathfrak l}\cap{\mathfrak k}$ of the Killing form of ${\mathfrak g}$, we have:
\begin{ex}\mbox{}\vspace*{-0.7cm}\\
\begin{itemize}
\item[(1)] $(G_1\times G_1,\Delta(G_1\times G_1), G_1\times\{e\})$: $\imath(\Omega_G)=2\Omega_L$.
\item[(2)] $(G_1\times G_1,\Delta(G_1\times G_1), G_1\times K_1)$: $\imath(\Omega_G)=2\Omega_L-\Omega_{L\cap K}$.
\item[(3)] $(SO_0(2,2n),SO_0(1,2n),U(1,n))$:  $\imath(\Omega_G)=2\Omega_L-\Omega_{L\cap K}$.
\item[(4)] $(SO_0(4,3),SO_0(4,1)\times SO(2),G_{2(2)})$:  $\imath(\Omega_G)=3\Omega_L-\frac{3}{2}\Omega_{L\cap K}+2\Omega_{\fl\cap\fs\cap\fq}$
\end{itemize}
\end{ex}
\noindent
In case (3), Schlichtkrull, Trapa and Vogan obtained a similar formula for $O(2n)/U(n)$ \cite{STV18}.
\begin{rem}
By inspection, we observe that $\imath(\Omega_G)$ involves `non-compact terms' terms from ${\mathfrak l}\cap{\mathfrak q}\cap{\mathfrak s}$ only when $(G,H,L)$ is not spherical, except in the first group case (1).
\end{rem}
%

\section{On the spectrum $\text{Spec}_{L^2(\Gamma\backslash G/H)}(D(G/H))$.}\label{spectrum}

Fix a transitive triple $(G,H,L)$. We have the following sequence of isomorphisms: 
\vspace*{0.5cm}
\begin{equation*}
L^2(\Gamma\backslash G/H)\simeq L^2(\Gamma\backslash L/L\cap H)\simeq L^2(\Gamma\backslash L)^{L\cap H}\simeq\widehat{\bigoplus}_{\pi\in\widehat{L}}V_{\widetilde{\pi},-\infty}^\Gamma\otimes V_{\pi}^{L\cap H}\hspace*{3cm}\text{ \small (Hilbert sum)}
\end{equation*}
\hspace*{3.8cm}$\bigcup$\hspace*{2.5cm}$\bigcup$\hspace*{3.2cm}$\bigcup$
%
\begin{equation*}
C^\infty(\Gamma\backslash G/H)\simeq C^\infty(\Gamma\backslash L/L\cap H)\simeq C^\infty(\Gamma\backslash L)^{L\cap H}\simeq\widehat{\bigoplus}_{\pi\in\widehat{L}}V_{\widetilde{\pi},-\infty}^\Gamma\otimes V_{\pi,\infty}^{L\cap H}\;\;\;\;\text{ \small (closure of algebraic direct sum)}
\end{equation*}
\hspace*{3.9cm}$\circlearrowleft$\hspace*{6.7cm}$\circlearrowleft$
\vspace*{-0.1cm}\\
\hspace*{3.5cm}$D\in D(G/H)$\hspace*{4.4cm}$\imath(D)\in D(L/L\cap H)$
\vspace*{0.2cm}\\
\noindent
Suppose now that $(G,H,L)$ is spherical, i.e $\dim\big(\text{Ind}_{P_L}^L \tau_L\otimes e^{\nu_L}\otimes 1\big)^{L\cap H}_{-\infty}<\infty$, for any $\tau_L\in\widehat{M_L}$ a finite dimensional irreducible representation of $M_L$ and $\nu_L$ a linear form on the complexified Lie algebra of $A_L$. Then by Proposition \ref{inv} and Casselman embedding Theorem, we deduce that
\begin{pro}\label{cass}
For any irreducible unitary representation $(\pi,V_\pi)$ of $L$, one has
\begin{equation*}
\dim V_{\pi,\infty}^{L\cap H}<\infty
\end{equation*}
\end{pro}
This means that the spectrum of $D(G/H)$ on $C^\infty(\Gamma\backslash G/H)$ may be reduced to the diagonalization of the operator $\pi(\imath(D))$ on the (infinitely many) finite dimensional blocks $V_{\pi,\infty}^{L\cap H}$. Then we obtain the following description of the spectrum pf $D(G/H)$ on $L^2(\Gamma\backslash G/H)$.
\begin{thm}\label{specthm} Suppose $(G,H,L)$ is a spherical triple.
\begin{itemize}
\item[(1)] There is a direct sum decomposition of $C^\infty(\Gamma\backslash G/H)$ and $L^2(\Gamma\backslash G/H)$ into joint eigenspaces of $D(G/H)$. $\text{Spec}_{L^2(\Gamma\backslash G/H)}(D(G/H))$ consists of the corresponding eigencharacters and their possible accumulation points.
\item[(2)] In the Lorentzian case $(SO_0(2,2n),SO_0(1,2n),U(1,n))$, $D(G/H)$ is generated by the Casimir $\Omega_G$, the set $\{\pi\in\widehat{L}\mid V_{\pi,\infty}^{L\cap H}\neq\{0\}\}$ is computed 
explicitly and the representations $\pi$ are identified as well as their contributions to eigenvalues for $\imath(\Omega_G)$. More precisely, if $\Delta$ denotes the Laplace operator on $\Gamma\backslash G/H$ induced by $\Omega_G$, one has
\begin{itemize}
\item $\spec(\Delta)\cap (-\infty,-n^2]$ comes from unitary principal series and, at $-n^2$, also from limits of discrete series of $L$. 
\item $\spec(\Delta)\cap (-n^2,0]$ consists of contributions from complementary series, ends of complementary series and non-integrable discrete series.
\item $\spec(\Delta)\cap (0,\infty)=\displaystyle\bigcup_{\ell=n+1}^\infty \{\ell^2-n^2\}$ is the contribution of integrable discrete series. The corresponding
eigenspaces are infinite dimensional.
\end{itemize}
\end{itemize}
\end{thm}
The last assertion about the contribution of integrable discrete series with $L\cap H$-invariants generalizes to arbitrary spherical triples $(G,H,L)$.
These contributions combine to infinite dimensional eigenspaces of $D(G/H)$ \cite{MOl}. The shortest proof of this uses Theorem \ref{ladmiss} instead of the embedding
$$
\imath:D(G/H)\hookrightarrow D(L/L\cap H)\ .$$ In view of Theorem \ref{ladmiss} the result can be rephrased in terms of $G$-representations as follows. Each integrable discrete series for $G/H$ contributes
an infinite dimensional eigenspace of $D(G/H)$ in $L^2(\Gamma\backslash G/H)$. Compare related results obtained in \cite{KK16}.

\section{Generalized matrix coefficients, unitarity and continuous spectrum.}\label{unitarity}

The Hilbert spaces $L^2(G/H)$ and $L^2(\Gamma\backslash G)$ are both unitary $G$-representations which decompose into a direct integral of irreducible ones. Though the Hilbert space $L^2(\Gamma\backslash G/H)$ is not a $G$-module, there is a natural way to involve representations of $G$. 

For a an admissible representation $(\rho,W_\rho)$ of $G$ of finite length, one may consider generalized matrix coefficients (see \cite{He14}): 
\begin{equation*}
W_{\widetilde{\rho},-\infty}\otimes W_{\rho,-\infty}\rightarrow C^{-\infty}(G),\;\widetilde{w}\otimes w \mapsto c_{\widetilde{w},w}
\end{equation*}
with 
\begin{eqnarray*}
c_{\widetilde{w},w}(f)&=&<\widetilde{w},\rho(f)w>\\
<\rho(f)w,v>&=&\int_G f(g)<\rho(g)w,v>dg\;\;\;\;\;\text{for }f\in C_c^\infty(G),\; v\in W_{\rho,\infty}
\end{eqnarray*}
Then there is a $D(G/H)$-equivariant map (which is injective if $\rho$ is irreducible):
\begin{equation}\label{gencoef}
W_{\widetilde{\rho},-\infty}^\Gamma\otimes W_{\rho,-\infty}^H\longrightarrow C^{-\infty}(\Gamma\backslash G/H)\simeq\widehat{\bigoplus}_{\pi\in\widehat{L}}V_{\widetilde{\pi},-\infty}^\Gamma\otimes V_{\pi,-\infty}^{L\cap H}
\end{equation}
\hspace*{5.3cm}$\circlearrowleft$\hspace*{6.5cm}$\circlearrowleft$
\vspace*{-0.1cm}\\
\hspace*{5cm}$D(G/H)$\hspace*{5cm}$D(L/L\cap H)$
\vspace*{0.2cm}\\

A comparison of (\ref{gencoef}) with Theorem \ref{specthm} leads, among other things, to the following funny observation. On one hand, from Theorem \ref{specthm}, there are complementary series representations $\pi$ of $U(1,n)$ contributing to finite dimensional eigenspaces. On the other hand, if $\rho$ were  unitary then $\rho_{\mid_L}$ would contain, by Theorem~\ref{ladmiss}, infinitely many unitary summands having $L\cap H$-invariants. Moreover, a result of Bergeron and Clozel \cite{BC}, Thm.~6.5.1, tell us that there are indeed lattices $\Gamma\subset U(1,n)$ such that $V_{\widetilde{\pi},-\infty}^\Gamma\ne\{0\}$ for complementary series $\pi$
with certain integral parameters. Therefore there must exist non-unitary representations $\rho$ of $G$ that are involved in $L^2(\Gamma\backslash G/H)$ via (\ref{gencoef}).

One can also try to use (\ref{gencoef}) to produce families of eigendistributions on $\Gamma\backslash G/H$ depending continuously on the spectral parameter. Let $0\ne \tilde v\in V_{\widetilde{\pi},-\infty}^\Gamma$ and assume we have a continuous nonconstant family  $t\mapsto \rho_t$, $t\in (a,b)\subset\mathbb R $ of  $H$-spherical principal series
of $G$ with $0\ne w_t\in W_{\rho_t,-\infty}^H$ and, the crucial ingredient,
a continuous family of non-zero intertwining operators $\Psi_t\in \Hom_L(W_{\rho_t,\infty},V_{\pi,\infty})$. Then $\Psi_t^*v\in  W_{\widetilde{\rho_t},-\infty}^\Gamma$ and
$$ t\mapsto c_{\Psi_t^*v,w_t} $$
is the desired continuous family of eigendistributions. Note that the existence of such a family for $\Gamma\backslash G/H$ excludes the existence of a discrete spectral decomposition as in
Theorem \ref{specthm}. Therefore such families can only exist for non-spherical triples $(G,H,L)$. One expects that such families are the building blocks for the continuous part
of the spectral decomposition of $L^2(\Gamma\backslash G/H)$. We are able to find such families of intertwiners for some non-spherical triples, and expect their existence
for all non-spherical $(G,H,L)$ except for the group case as in Example \ref{obst} (1). For instance, for the non-spherical triple $(SO(8,\mathbb C),SO_0(1,7), Spin(7,\mathbb C))$
the existence of many families $\Psi_t$ is ensured by the work of Jan Frahm \cite{Mo17}.


\end{document}